\documentclass[a4paper,10pt]{article}
\usepackage{amsmath}
\usepackage{amssymb}
\usepackage{amsthm}
\theoremstyle{plain}
\numberwithin{equation}{section}
\newtheorem{theo}{Theorem}[section]
\newtheorem{cor}{Corollary}[section]
\newtheorem{lem}{Lemma}[section]
\theoremstyle{definition}
\newtheorem{defi}{Definition}[section]
\theoremstyle{remark}
\newtheorem{rem}{Remark}[section]

\begin{document}
\begin{Large}
 {\bf Subexponential densities of compound Poisson sums and the supremum of a random walk}
\end{Large}
\begin{center}
Takaaki Shimura\footnote[1]{Takaaki Shimura \\
The Institute of Statistical Mathematics, 10-3 Midori-cho, Tachikawa, Tokyo 190-8562, Japan\\
E-mail: shimura@ism.ac.jp \\} $\quad$ Toshiro Watanabe \footnote[2]{
Toshiro Watanabe\\
Center for Mathematical Sciences, The University of Aizu, Fukushima 965-8580, Japan \\
E-mail: t-watanb@u-aizu.ac.jp}
\end{center}  
\medskip
{\bf Abstract} 
We characterize the subexponential densities on $(0,\infty)$ for compound Poisson
    distributions on $[0,\infty)$ with absolutely continuous L\'evy measures. As a corollary, we show that the class of all subexponential probability density functions on $\mathbb R_+$ is closed under generalized convolution roots of compound Poisson sums. Moreover, we give an   application to the subexponential density on $(0,\infty)$ for the distribution of the supremum of a random walk.
  
\medskip

  compound Poisson distribution, random walk\\
 

\section{Introduction and main results}
\medskip

  In what follows, we denote by $\mathbb R$ the real line and by $\mathbb R_{+}$ the half line $[0,\infty)$. Denote by $\mathbb N$ the totality of positive integers. The symbol $\delta_a(dx)$ stands for the delta measure at $a \in \mathbb R$. Let $\eta$ and $\rho$ be probability distributions on $\mathbb R$. We denote by $\eta*\rho$ the convolution of $\eta$ and $\rho$ and by $\rho^{n*}$
$n$-th convolution power of $\rho$ with the understanding that $\rho^{0*}(dx)=\delta_0(dx)$. 
Let $f(x)$ and $g(x)$ be probability density functions on $\mathbb R$. We denote by $f \otimes g(x)$ the convolution of $f(x)$ and $g(x)$ and by $f^{n\otimes}(x)$ $n$-th convolution power of $f(x)$ for $n \in \mathbb N$. For positive functions $f_1(x)$ and $g_1(x)$ on $[A,\infty)$ for some $A \in \mathbb R$, we define the relation $f_1(x)\sim g_1(x)$ by $\lim_{x \to \infty}f_1(x)/g_1(x)=1.$  We use the symbols $\mathcal{L}$ and $\mathcal{S}$ in the sense of long-tailed and subexponential, respectively.
\begin{defi}
(i) A nonnegative measurable function $g(x)$ on $\mathbb R$ belongs to the class ${\bf L}$ if $g(x+a)\sim g(x)$ for every $ a \in \mathbb R$.

(ii) A probability density function $g(x)$ on $\mathbb R$ belongs to the class $\mathcal{L}_d$ if $g(x) \in {\bf L}$.

(iii)  A probability density function $g(x)$ on $\mathbb R$ belongs to the class $\mathcal{S}_d$ if $g(x) \in \mathcal{L}_d$ and $g^{2\otimes}(x)\sim 2g(x)$.  
\end{defi}
\begin{defi} (i)  Let $\Delta:=(0,c]$ with $c >0$. A distribution $\rho$ on $\mathbb R$ belongs to the class $\mathcal{L}_{\Delta}$ if $\rho((x,x+c]) \in {\bf L}$. A distribution $\rho$ on $\mathbb R$ belongs to the class $\mathcal{L}_{loc}$ if $\rho \in \mathcal{L}_{\Delta}$ for each $\Delta:=(0,c]$ with $c >0$.

(ii)  Let $\Delta:=(0,c]$ with $c >0$. A distribution $\rho$ on $\mathbb R$ belongs to the class $\mathcal{S}_{\Delta}$ if $\rho \in \mathcal{L}_{\Delta}$ and  $\rho^{2*}((x,x+c])\sim 2\rho((x,x+c])$. 
 A distribution $\rho$ on $\mathbb R$ belongs to the class $\mathcal{S}_{loc}$ if $\rho \in \mathcal{S}_{\Delta}$ for each $\Delta:=(0,c]$ with $c >0$.

\end{defi}
\medskip
  Functions in the class ${\bf L}$ are called {\it long-tailed functions}. Probability density functions in the classes  $\mathcal{L}_d$ and  $\mathcal{S}_d$ are called {\it long-tailed densities} and {\it subexponential densities}, respectively. Density functions, which are not necessarily probability ones, are called subexponential if their normalized probability densities are subexponential. Note that if $f(x) \in \mathcal{L}_d$, then $\lim_{x \to \infty}f(x)=0$ and $\lim_{x \to \infty}e^{sx}f(x)=\infty$ for every $s >0$. See Foss et al. (2013). Distributions in the classes  $\mathcal{S}_{\Delta}$ and $\mathcal{S}_{loc}$ are called $\Delta$-{\it subexponential} and {\it locally subexponential}, respectively. The class $\mathcal{S}_{\Delta}$ was introduced by  Asmussen et al. (2003).

 Let $\mu$ be a compound Poisson distribution on $\mathbb R_+$ with  L\'evy measure $\nu$.  Denote by $\mu^{t*}$ $t$-th convolution power of $\mu$ for $t >0$. Compound Poisson distributions have many applications. Later, we give an   application to the subexponential density on $(0,\infty)$ for the distribution  $\pi$ of the supremum of a random walk. The compound Poisson distribution $\pi$ is important in classical ruin theory and queueing theory.   It is obvious that a compound Poisson distribution $\mu$ on $\mathbb R_+$ is absolutely continuous on $(0,\infty)$ if and only if its L\'evy measure $\nu$ is absolutely continuous. In the following theorem and the corollary, let $\mu$ be a compound Poisson distribution on $\mathbb R_+$ with absolutely continuous L\'evy measure $\nu$. Let $\nu(dx) = \lambda\phi(x)dx$ on $\mathbb R_+$ with 
$\lambda:=\nu((0,\infty)) \in (0,\infty)$ and define a probability density function of compound Poisson sums $p^t(x)$ on $\mathbb R_+$ for $t >0$ as
\begin{equation}
p^t(x):=(e^{\lambda t}-1)^{-1}\sum_{n=1}^{\infty}\frac{(\lambda t)^n}{n!}\phi^{n\otimes}(x) \nonumber
\end{equation}
and let $p(x):=p^1(x)$. Then, we have
\begin{equation}
\mu^{t*}(dx)=e^{-\lambda t}\delta_0(dx)+(1-e^{-\lambda t})p^t(x)dx. \nonumber
\end{equation}
We are concerned with the asymptotic relation between the densities  $p(x)$ and $\phi(x)$. Namely, we consider the following problem.
\medskip
\newline
{\bf Problem 1}  Are the following assertions equivalent ?

(a) $p(x)  \in \mathcal{S}_d$.

(b)  $\phi(x)  \in \mathcal{S}_d$.

(c)  $\phi(x)  \in \mathcal{L}_d$
 and
\begin{equation}
(1-e^{-\lambda })p(x)\sim \lambda\phi(x). \nonumber
\end{equation}

(d) $\phi(x)  \in \mathcal{L}_d$
 and there is $C \in (0,\infty)$ such that
\begin{equation}
p(x)\sim C\phi(x). \nonumber
\end{equation}
\medskip
\newline
 We say that the class $\mathcal{S}_d$ of probability density functions on $\mathbb R_+$ is closed under generalized convolution roots of compound Poisson sums if assertion (a) implies (b). We answer Problem 1 under the assumption that  $\int_0^{\infty}(\phi(x))^2dx < \infty.$
 \medskip
\begin{theo} Let assertions (a)-(d) be  the same as those in  Problem 1. Then, we have the following. 

(i) Assertion (a) implies (b).

(ii) Assertions (c) and (d) are equivalent and (c) implies (a).

(iii) Assume that $\int_0^{\infty}(\phi(x))^2dx < \infty.$ Then, (b) implies (c). Thus, under the assumption that $\int_0^{\infty}(\phi(x))^2dx < \infty,$ all assertions (a)-(d) are equivalent.
\end{theo}
\medskip
\begin{rem} In the theorem above,  without the assumption that $\int_0^{\infty}(\phi(x))^2dx < \infty$, assertion (b) does not necessarily imply (c). See Remark 2.1 (b) below.   Kl\"uppelberg (1989) also obtained in her Corollary 3.3  a result analogous to Theorem 1.1 under the assumptions that $\phi(x) \in \mathcal{L}_d$ and that $\phi(x)$ is bounded on $\mathbb R_+$. We do not need the assumption  that $\phi(x) \in \mathcal{L}_d$ for the proof of (i) of Theorem 1.1.  Thus, the class $\mathcal{S}_d$ of probability density functions on $\mathbb R_+$ is closed under generalized convolution roots of compound Poisson sums. As an example, 
every compound Poisson distribution on $\mathbb R_+$ has a square integrable regularly varying density on $(0,\infty)$
 if and only if so does $\nu$.
\end{rem}
We can characterize the precise asymptotic behavior of the distributions of compound Poisson processes with absolutely continuous L\'evy measures as follows.
\begin{cor}  Assume that $\int_0^{\infty}(\phi(x))^2dx < \infty.$ Then, we have the following. 

(i) If $p^t(x) \in \mathcal{S}_d$ for some $t >0$, then $p^t(x) \in \mathcal{S}_d$ for all $t >0$ and
$$(1-e^{-\lambda t }) p^t(x) \sim t(1-e^{-\lambda }) p(x). $$

(ii) If $p(x) \in \mathcal{L}_d$ and, for some 
$t \in  (0,1) \cup (1,\infty)$, there is 
$ C(t) \in (0,\infty)$            
 such that
\begin{equation}
(1-e^{-\lambda t }) p^t(x) \sim C(t)(1-e^{-\lambda }) p(x), 
\end{equation}
then $C(t)=t$ and $p(x) \in \mathcal{S}_d$.
\end{cor}
\medskip

The organization of this paper is as follows. In Sect.\ 2, we explain basic results on the classes  $\mathcal{S}_d$
and $\mathcal{S}_{loc}$ as preliminaries. In Sect.\ 3, we prove Theorem 1.1 together with its corollary. In Sect.\ 4, we give an application of Theorem 1.1 to the supremum of a random walk.

\section{Preliminaries}

 In this section, we give several fundamental results on the classes  $\mathcal{S}_d$
and $\mathcal{S}_{loc}$.

\begin{lem} Let $f(x)$ and $g(x)$ be probability density functions on $\mathbb R$.

(i) If $f(x),g(x) \in \mathcal{L}_d$, then $f\otimes g(x) \in \mathcal{L}_d$.

(ii) Let $f(x) \in \mathcal{L}_d$ and define a distribution $\rho$ on $\mathbb R$ by 
$$\rho(dx):=c_0\delta_0(dx) +(1-c_0)f(x)dx$$ 
with $c_0 \in [0,1).$ Then, $\rho \in \mathcal{S}_{loc}$ if and only if $f(x) \in \mathcal{S}_d$. 
\end{lem}
\medskip

Proof.   Assertion (i) is due to Theorem 4.3 of Foss et al. (2013). Next, we prove assertion (ii). Assume that $f(x) \in \mathcal{L}_d$. Then, by (i), $f^{2\otimes}(x) \in \mathcal{L}_d$
and hence $f(x+u) \sim f(x)$ and $f^{2\otimes}(x+u) \sim f^{2\otimes}(x)$ uniformly in $u \in [0,c]$ with $c >0$. Thus, we have, for $x > 0$,
\begin{equation}
\begin{split}
&\rho^{2*}((x,x+c])\\ \nonumber
&=2c_0(1-c_0)\int_0^{c}f(x+u)du+(1-c_0)^2\int_0^{c}f^{2\otimes}(x+u)du\\
&\sim 2c_0(1-c_0)cf(x)+(1-c_0)^2cf^{2\otimes}(x).
\end{split}
\end{equation}
Hence, we see that  $\rho \in \mathcal{S}_{loc}$ if and only if $f^{2\otimes}(x) \sim 2f(x)$, namely, $f(x) \in \mathcal{S}_d$. 
$\Box$
\begin{lem}  Let $f(x)$ and $g(x)$ be probability density functions on $\mathbb R_+$.

(i) If $f(x) \in \mathcal{S}_d$ and $g(x) \sim cf(x)$ with $c \in (0,\infty)$, then $g(x) \in \mathcal{S}_d$.

(ii) Assume that $f(x) \in \mathcal{S}_d$ and $\int_0^{\infty}(f(x))^2dx < \infty$. Then, for any $\epsilon > 0$, there are $x_0(\epsilon) >0$ and $C(\epsilon) >0$ such that, for all $x > x_0(\epsilon)$ and all $n \in \mathbb N$, 
$$f^{n\otimes}(x) \leq C(\epsilon)(1+\epsilon)^nf(x).$$

(iii) If $f(x) \in \mathcal{S}_d$, then, for all $n \in \mathbb N$,
$$ f^{n\otimes}(x) \sim nf(x).$$
\end{lem}

\medskip
Proof. Assertions (i) and (iii) are due to Theorem 4.8 and Corollary 4.10 of Foss et al. (2013), respectively. Assertion (ii) is a modification of Theorem 4.11 of Foss et al. (2013). \hfill $\Box$

\medskip
\begin{rem} (a)  Watanabe and Yamamuro (2017) showed in Theorem 1.2  that assertion (i) is not necessarily true for probability density functions $f(x)$ and $g(x)$ on $\mathbb R$.
Assertion (ii) is called Kesten's bound. In Theorem 4.11 of Foss et al. (2013), boundedness of $f(x)$ on $\mathbb R_+$ is assumed in place of square integrability of $f(x)$ on $\mathbb R_+$.
 If $C:=\int_0^{\infty}(f(x))^2dx < \infty$, then $f^{n\otimes}(x)\leq C$ for all $x >0$ and all integers $n \geq 2$. Thus, without any change, the proof of Theorem 4.11 of Foss et al. (2013) leads to assertion (ii).  See Finkelshtein and Tkachov (2018)  
  for Kesten's bound for a probability density function $f(x)$ on $\mathbb R$. It needs some additional conditions for $f(x)$  on $\mathbb R$.
Assertion (iii) goes back to Chover et al. (1973). 

(b) We show  that there exists a subexponential probability density function $\phi(x)$ on $\mathbb R_+$ which does not satisfy Kesten's bound.  Let 
$$h(x):=1_{(-e^{-1},e^{-1})}(x)|x|^{-1}|\log |x||^{-2}$$
 and let $g(x):=2^{-1}h(x-1)$. Then $g(x)$ is a probability density function on $\mathbb R_+$.  We have easily, for all $n  \in \mathbb N$, 
$$\lim_{x \to n}g^{n \otimes}(x)=\infty.$$
Let $f(x) \in \mathcal{S}_d$ be bounded  on $\mathbb R_+$ and let $\phi(x):=2^{-1}f(x)+2^{-1}g(x)$. Then we have $\phi(x) \in \mathcal{S}_d$ but Kesten's bound does not holds for $\phi(x)$. Moreover, we see that, for all $n  \in \mathbb N$, 
$$\lim_{x \to n}p(x)=\lim_{x \to n}(e-1)^{-1}\sum_{k=1}^{\infty}\frac{\phi^{k \otimes}(x)}{k!}=\infty.$$
Thus, assertion (b) does not necessarily  implies (c) in Problem 1 without the assumption that $\int_0^{\infty}(\phi(x))^2dx < \infty.$
\end{rem}

Watanabe and Yamamuro (2010) used the principal results of Watanabe (2008) on the convolution equivalence of infinitely divisible distributions on $\mathbb R$  to prove the following  lemmas. Our main results essentially depend on those two  results.
\begin{lem} (Theorem 1.1 of Watanabe and Yamamuro (2010)) Let $\mu$ be an  infinitely divisible distribution on $\mathbb R_+$  with L\'evy measure $\nu$. Then, the following are equivalent : 

(1) $\mu \in \mathcal{S}_{loc}$.

(2) $\nu_{(1)} \in \mathcal{S}_{loc}$.

(3) $\nu_{(1)} \in \mathcal{L}_{loc}$
and $\mu((x,x+c]) \sim \nu((x,x+c])$ for all $c >0$.

(4) $\nu_{(1)} \in \mathcal{L}_{loc}$
and there is $C \in (0,\infty)$ such that $\mu((x,x+c]) \sim C\nu((x,x+c])$ for all $c >0$.
\end{lem}

\begin{lem} (Theorem 1.2 of Watanabe and Yamamuro (2010)) Let $\mu$ be an  infinitely divisible distribution on $\mathbb R_+$  with L\'evy measure $\nu$. Then, we have the following.

(i) If $\mu^{t*} \in \mathcal{S}_{loc}$ for some $t >0$, then $\mu^{t*} \in \mathcal{S}_{loc}$ for all $t >0$ and
$$ \mu^{t*}((x,x+c]) \sim t\mu((x,x+c]) $$
 for all $t >0 $ and for all $c >0. $

(ii) If $\mu \in \mathcal{S}_{loc}$ and, for some 
$t \in  (0,1) \cup (1,\infty)$, there is 
$ C(t) \in (0,\infty)$            
 such that

\begin{equation}
 \mu^{t*}((x,x+c]) \sim C(t)\mu((x,x+c])
\end{equation}
 for all $c >0, $ then $C(t)=t$ and $\mu \in \mathcal{S}_{loc}$.
\end{lem}

 \section{
 Proofs of Theorem 1.1 and its corollary}
 
 \medskip
 
 For an integrable function $h(x)$ on $\mathbb R_+$, denote by $L_h(t)$ the Laplace transform of $h(x)$, namely, for  $t \in \mathbb R_+$,
 $$L_h(t):=\int_0^{\infty}e^{-tx}h(x)dx.$$ 
We begin with a key lemma for the proof of Theorem 1.1.

\medskip

\begin{lem}  Let $\mu$ be a compound Poisson distribution with absolutely continuous l\'evy measure $\nu$. Let $\nu(dx)=\lambda\phi(x)dx$ on $\mathbb R_+$ with $\lambda:=\nu((0,\infty)) \in (0,\infty)$ and let $\lambda_1\phi_1(x):=\lambda 1_{(c_1,\infty)}(x)\phi(x)$ with $c_1 >0$ and $\lambda_1:=\lambda\int_{c_1}^{\infty}\phi(x)dx.$ Define probability density functions $p(x)$ and $p_1(x)$ on $\mathbb R_+$ as
$$p(x):=(e^{\lambda }-1)^{-1}\sum_{n=1}^{\infty}\frac{\lambda^n}{n!}\phi^{n\otimes}(x)$$
and
$$p_1(x):=(e^{\lambda_1 }-1)^{-1}\sum_{n=1}^{\infty}\frac{\lambda_1^n}{n!}\phi_1^{n\otimes}(x).$$
Then, $p(x) \in \mathcal{S}_d$ implies that $p_1(x) \in \mathcal{S}_d$ for every $c_1>0$.
\end{lem}
\medskip

Proof. Suppose that $p(x) \in \mathcal{S}_d$. If $\nu((0,c_1])=0,$ then $p_1(x)=p(x) \in \mathcal{S}_d$. Hence, we can assume that 
$\nu((0,c_1])>0.$ We define a strictly increasing finite sequence
 $\{a_n\}_{n=0}^{N}$ with $N \in \mathbb N$ such that $a_0=0$ and $a_N=c_1$. 
 Let
 $$\alpha_n:=\exp(-\nu((a_n,\infty))) \in (0,1)$$
 for $0 \leq n \leq N$ and let
  $$\beta_n:=\exp(-\nu((a_{n-1},a_n])) $$
for $1 \leq n \leq N$. Then, we have $\alpha_n \beta_n =\alpha_{n-1}$
for $1 \leq n \leq N$. We can choose $\beta_n$ such that, for $1 \leq n \leq N$,
\begin{equation}
2^{-1} < \beta_n < 1.
\end{equation}
Define probability density functions $\varphi_n(x)$ and $f_n(x)$ on $\mathbb R_+$ as, for $0 \leq n \leq N$,
$$(-\log \alpha_n)\varphi_n(x):=\lambda 1_{(a_n,\infty)}(x)\phi(x),$$
and
$$(1-\alpha_n)f_n(x):=\alpha_n\sum_{k=1}^{\infty}\frac{(-\log \alpha_n)^k}{k!}\varphi_n^{k \otimes}(x).$$
Moreover, define probability density functions $\psi_n(x)$ and $g_n(x)$ on $\mathbb R_+$ as, for $1 \leq n \leq N$,
$$(-\log \beta_n)\psi_n(x):=\lambda 1_{(a_{n-1},a_n]}(x)\phi(x),$$
and
\begin{equation}
(1-\beta_n)g_n(x):=\beta_n\sum_{k=1}^{\infty}\frac{(-\log \beta_n)^k}{k!}\psi_n^{k \otimes}(x).
\end{equation}
 Since, for $1 \leq n \leq N$,
 $$(-\log \alpha_{n-1})\varphi_{n-1}(x)= (-\log \beta_n)\psi_n(x) + (-\log \alpha_n)\varphi_n(x),$$
 we have, for $x \in \mathbb R_+$,
\begin{equation}
\begin{split}
 (1-\alpha_{n-1})f_{n-1}(x)=&\alpha_n(1-\beta_n)g_n(x)+\beta_n(1-\alpha_n)f_n(x)\\
& +(1-\alpha_n)(1-\beta_n)f_n\otimes g_n(x).
\end{split}
\end{equation}
We shall prove that if $f_{n-1}(x) \in \mathcal{S}_d$ for some $1 \leq n \leq N$, then $f_n(x) \in \mathcal{S}_d$. Suppose that $f_{n-1}(x) \in \mathcal{S}_d$ for some $1 \leq n \leq N$.
Define constants $C^*$ and $C_*$ as
$$C^*:=\limsup_{x \to \infty}\frac{f_n(x)}{f_{n-1}(x)},\quad  C_*:=\liminf_{x \to \infty}\frac{f_n(x)}{f_{n-1}(x)}.$$
We find from (3.3) that
$$ 0 \leq C_*\leq C^*< \infty.$$
By virtue of Fatou's lemma, we have
\begin{equation}
\begin{split}
&\liminf_{x \to \infty}\frac{f_n\otimes g_n(x)}{f_{n-1}(x)}\\ \nonumber
&\geq\lim_{M \to \infty}\int_0^M\liminf_{x \to \infty}\frac{f_n(x-u)}{f_{n-1}(x-u)}\frac{f_{n-1}(x-u)}{f_{n-1}(x)}g_n(u)du\\
&\geq  C_*\int_0^{\infty}g_n(u)du=  C_*.
\end{split}
\end{equation}
Thus, we obtain from (3.3) that
\begin{equation}
1-\alpha_{n-1}\geq \beta_n(1-\alpha_n)C^*+(1-\alpha_n)(1-\beta_n)C_*.
\end{equation}
Define a probability density function $h_n(x)$
 on $\mathbb R_+$ for $0\leq n \leq N$ as
$$(1-\alpha_{n}^2)h_n(x):=\alpha_n^2\sum_{k=1}^{\infty}\frac{(-2\log \alpha_n)^k}{k!}\varphi_n^{k \otimes}(x).$$
Then, we have for $1 \leq n \leq N$
\begin{equation}
(1-\alpha_{n-1}^2)h_{n-1}(x)\geq \alpha_{n-1}^2\sum_{k=1}^{\infty}\frac{(-2\log \beta_n)^k}{k!}\psi_n^{k \otimes}(x).
\end{equation}
Let $M >0$ and $K(M):=[M/a_n].$   Here, the symbol $[x]$ stands for the largest integer not exceeding $x \in \mathbb R$.   Since $\psi_n^{k \otimes}(x)=0$ for $1 \leq k \leq K(M)$ and $x > M$, we see from (3.2) and (3.5) that, for $ x > M$, 
\begin{equation}
g_n(x)\leq \frac{(1-\alpha_{n-1}^2)\beta_n}{\alpha_{n-1}^2(1-\beta_n)
2^{K(M)}}h_{n-1}(x).
\end{equation}
Since we assume that $f_{n-1}(x) \in \mathcal{S}_d$, note that
\begin{equation}
\begin{split}
(1-\alpha_{n-1}^2)h_{n-1}(x)&= 2\alpha_{n-1}(1-\alpha_{n-1})f_{n-1}(x)+ (1-\alpha_{n-1})^2f_{n-1}^{2\otimes}(x)\\ \nonumber
&\sim 2(1-\alpha_{n-1})f_{n-1}(x).
\end{split}
\end{equation}
Since $\lim_{M \to \infty}K(M)= \infty$, we obtain from (3.6) that
 \begin{equation}
 \lim_{x \to \infty}\frac{g_n(x)}{f_{n-1}(x)}=0.
 \end{equation}
 Let
 $$f_n\otimes g_n(x) = I_1(x) +I_2(x),$$
 where
 $$I_1(x):=\int_0^Mf_n(x-u)g_n(u)du$$
 and
 $$I_2(x):=\int_M^xf_n(x-u)g_n(u)du.$$
Then, we have by Fatou's lemma
\begin{equation}
\begin{split}
&\limsup_{x \to \infty}\frac{I_1(x)}{f_{n-1}(x)}\\
 & \leq \int_0^M\limsup_{x \to \infty}\frac{f_n(x-u)}{f_{n-1}(x-u)}\frac{f_{n-1}(x-u)}{f_{n-1}(x)}g_n(u)du\\
& =  C^*\int_0^{M}g_n(u)du.
\end{split}
\end{equation}
We find from (3.3) that, for $x >0$,
$$f_n(x) \leq \frac{1-\alpha_{n-1}}{\beta_n (1-\alpha_n)}f_{n-1}(x).$$
Note from (3.7) that there is $\epsilon (M)>0$ such that
$\lim_{M \to\infty}\epsilon (M)=0$ and $g_n(x)\leq \epsilon (M)f_{n-1}(x)$ for $x >M$. Thus, we have
\begin{equation}
\begin{split}
&\limsup_{x \to \infty}\frac{I_2(x)}{f_{n-1}(x)}\\ \nonumber
 & \leq \limsup_{x \to \infty}\frac{1}{f_{n-1}(x)}\int_M^x\frac{1-\alpha_{n-1}}{\beta_n (1-\alpha_n)}f_{n-1}(x-u)\epsilon (M)f_{n-1}(u)du\\
& \leq\frac{\epsilon (M) (1-\alpha_{n-1})}{\beta_n (1-\alpha_n)}\lim_{x \to\infty}\frac{f_{n-1}^{2\otimes}(x)}{f_{n-1}(x)}\\
&=\frac{2\epsilon (M) (1-\alpha_{n-1})}{\beta_n (1-\alpha_n)}.
\end{split}
\end{equation}
Hence, we see from (3.8) that
\begin{equation}
\begin{split}
&\limsup_{x \to \infty}\frac{f_n\otimes g_n(x)}{f_{n-1}(x)}\\
 & \leq \lim_{M \to \infty}\left(C^*\int_0^M g_n(u)du+ \frac{2\epsilon (M) (1-\alpha_{n-1})}{\beta_n (1-\alpha_n)}\right)=C^*.
\end{split}
\end{equation}
Thus, we obtain from (3.7) and (3.9) that
\begin{equation}
1-\alpha_{n-1}\leq \beta_n(1-\alpha_n)C_*+(1-\alpha_n)(1-\beta_n)C^*.
\end{equation}
Hence, $C^* >0$ and by (3.1),  (3.4), and (3.10) we have
$$0\geq (C^*-C_*)(1-\alpha_n)(2\beta_n-1)\geq 0.$$
Thus, we have $C^*=C_* \in (0,\infty)$ and $f_n(x)\sim C^*f_{n-1}(x).$ By  (i) of Lemma 2.2, we have proved that  $f_n(x) \in \mathcal{S}_d$. Since $f_0(x)=p(x) \in \mathcal{S}_d$, we conclude that $f_N(x)=p_1(x) \in \mathcal{S}_d$ by induction.\hfill $\Box$

\medskip

Proof of Theorem 1.1. First, we prove assertion (i). Suppose that (a) holds. Since $\lim_{x \to \infty}p(x) =0,$
we have $\lim_{x \to \infty}\phi(x) =0.$ Choose sufficiently large $c_1>0$ such that $e^{\lambda_1} <2$ with $\lambda_1:=\nu((c_1,\infty))$ and $\sup_{x>c_1}\phi(x) < \infty$. Then, we see from Lemma 3.1 that $ p_1(x) \in \mathcal{S}_d$ and $p_1(x)$ is bounded on $\mathbb R_+$. Noting that $0 < e^{\lambda_1}-1 <1$, define a function $\phi_0(x)$ on $\mathbb R_+$ as
$$ \lambda_1\phi_0(x):=-\sum_{n=1}^{\infty}\frac{(1-e^{\lambda_1})^n}{n}p_1^{n\otimes}(x).$$
Since $ p_1(x) \in \mathcal{S}_d$ and $p_1(x)$ is bounded on $\mathbb R_+$, we obtain from (ii) and (iii) of Lemma 2.2 that
\begin{equation}
\begin{split}
\lambda_1\phi_0(x)&\sim -\sum_{n=1}^{\infty}(1-e^{\lambda_1})^n p_1(x)\\
&=(1-e^{-\lambda_1}) p_1(x).
\end{split}
\end{equation}
Let $\tilde \phi_0(x):=\phi_0(x)\vee 0.$ By using Remark 21.6 of Sato (2013), we have, for $t \in \mathbb R_+$,
\begin{equation}
\begin{split}
\lambda_1L_{\phi_0}(t)&= -\sum_{n=1}^{\infty}\frac{(1-e^{\lambda_1})^n}{n} (L_{p_1}(t))^n\\ \nonumber
&=\log(1-(1-e^{\lambda_1})L_{ p_1}(t))\\
&=\log(e^{\lambda_1}\exp(\lambda_1\int_0^{\infty}(e^{-tx}-1)\phi_1(x)dx))\\
&=\lambda_1L_{\phi_1}(t).
\end{split}
\end{equation}
Thereby, we see that $\tilde \phi_0(x):=\phi_1(x)$ for a.e. $ x \in \mathbb R_+.$ Thus, we obtain from (3.11) and (i) of Lemma 2.2 that 
$ \tilde \phi_0(x) \in \mathcal{S}_d$ and $ \tilde \phi_0(x)$ is bounded on $\mathbb R_+$. Hence, by (ii) and (iii) of Lemma 2.2 and (3.11), we have
\begin{equation}
\begin{split}
\lambda_1\phi_1(x)&=(e^{\lambda_1}-1)p_1(x)-\sum_{n=2}^{\infty}\frac{\lambda_1^n}{n!}\tilde \phi_0^{n\otimes}(x)\\ \nonumber
& \sim (1-e^{-\lambda_1})p_1(x).  
\end{split}
\end{equation}
Since 
$$ \phi(x)\sim \frac{\lambda_1}{\lambda}\phi_1(x) \sim \frac{1-e^{-\lambda_1}}{\lambda}p_1(x),$$
we find from (i) of Lemma 2.2 that $ \phi(x) \in \mathcal{S}_d$. Next, we see from Lemma 2.3 and  (ii) of Lemma 2.1 that assertion (ii) is valid and from (ii) and (iii) of Lemma 2.2 that assertion (iii) is also valid.\hfill $\Box$

\medskip

Proof of Corollary 1.1. Assume that $\int_0^{\infty}(\phi(x))^2dx < \infty.$ First, we prove assertion (i). Let $\nu(dx)=\lambda\phi(x)dx$
 with $\lambda:=\nu((0,\infty))\in (0,\infty)$. Suppose that  $ p^t(x) \in \mathcal{S}_d$ for some $t >0$. Then, we obtain from Theorem 1.1 that  $ \phi(x) \in \mathcal{S}_d$ and hence  $ p^t(x) \in \mathcal{S}_d$ for all $t >0$. Moreover, we find again from Theorem 1.1 that, for all $t >0$,
 $$(1-e^{-\lambda t})p^t(x) \sim t\lambda\phi(x) \sim t(1-e^{-\lambda })p(x).$$
 Next, we prove assertion (ii). Suppose that  $ p(x) \in \mathcal{S}_d$ and, for some $t \in (0,1)\cup(1,\infty)$, there is $C(t) \in (0,\infty)$ such that (1.1) holds. Then, we have $\mu \in \mathcal{L}_{loc}$ and (2.1) holds. Thus, we obtain from Lemma 2.4 that $C(t)=t$ and $\mu \in \mathcal{S}_{loc}$ and hence, by $ p(x) \in \mathcal{L}_d$ and (ii) of Lemma 2.1, $ p(x) \in \mathcal{S}_d$.\hfill $\Box$
  
 \section{Application to the supremum of a random walk}
 
  In this section, we characterize the subexponentiality of the density on $(0,\infty)$ for the distribution of the supremum of a random walk. For a distribution $\rho$ on $\mathbb R$, denote the tail of $\rho$ by $\bar\rho(x):=\rho((x,\infty))$ for $x \in \mathbb R$.
 \begin{defi} (i) A   distribution $\rho$ on $\mathbb R$ belongs to the class $\mathcal{L}$ if $\bar\rho(x) \in \bf L$. A probability distribution $\rho$ on $\mathbb R$ is called {\it subexponential} if  $\rho\in \mathcal{L}$ and  
 $$\overline{\rho^{2*}}(x)\sim 2\bar\rho(x).$$ 
 The class of all subexponential 
 distributions on $\mathbb R$ is denoted by $\mathcal{S}$. 
 
 (ii)  A   distribution $\rho$ on $\mathbb R$ belongs to the class $\mathcal{S}^*$ if $\bar\rho(x) > 0$ for all  $x \in\mathbb R$, $\int_0^{\infty}x\rho(dx)<\infty$ and if 
 $$ \int_0^x\bar\rho(x-y)\bar\rho(y)dy \sim 2 \int_0^{\infty}u\rho(du)\bar\rho(x).$$
\end{defi}
\medskip

 Refer to Kl\"uppelberg (1988) for the class  $\mathcal{S}^*$.                    Note that the class  $\mathcal{S}^*$ is included in the class   $\mathcal{S}$. We define a function $\varphi_{\alpha,\lambda}(s)$ on $(-1/\lambda,1/\lambda)$ for $0 < \lambda < 1$ and $\alpha >0$ by
\begin{equation}
 \varphi_{\alpha,\lambda}(s):=\left( \frac{1-\lambda}{1-\lambda s}\right)^{\alpha}=\sum_{n=0}^{\infty}{\alpha+n-1\choose\alpha-1}(1-\lambda)^{\alpha}\lambda^ns^n. \nonumber
 \end{equation}
 Define by $\varphi_{\alpha,\lambda}'(s)$ the derivative of $\varphi_{\alpha,\lambda}(s)$. Let 
 $$c_0:=(1-\lambda)^{\alpha}\mbox{ and }\delta:=-\log(1-\lambda).$$
 
  \begin{lem}  Let $0 < \lambda < 1$ and $\alpha >0$ and let $f(x)$ be a probability density function on $\mathbb R_+$. Define probability density functions $p(x)$ and  $\phi(x)$ on $\mathbb R_+$ as
 \begin{equation}
 p(x):=\frac{c_0}{1-c_0}\sum_{n=1}^{\infty}{\alpha+n-1\choose\alpha-1}\lambda^nf^{n\otimes}(x)
 \end{equation}
and
 \begin{equation}
 \phi(x):=\delta^{-1}\sum_{n=1}^{\infty}\frac{\lambda^n}{n}f^{n\otimes}(x).
 \end{equation} 
Then, we see, for all $x \in \mathbb R_+$, 
\begin{equation}
p(x)=(e^{\alpha\delta}-1)^{-1}\sum_{n=1}^{\infty}\frac{(\alpha\delta)^n}{n!}\phi^{n\otimes}(x).
\end{equation} 
\end{lem}
\medskip

Proof. Let 
$$\tilde p(x):=(e^{\alpha\delta}-1)^{-1}\sum_{n=1}^{\infty}\frac{(\alpha\delta)^n}{n!}\phi^{n\otimes}(x).$$
Substituting (4.2) in the above equation, we have
$$ \tilde p(x) = \sum_{n=1}^{\infty}p_nf^{n\otimes}(x),$$
for some nonnegative sequence $\{p_n\}_{n=1}^{\infty}$ with $\sum_{n=1}^{\infty}p_n =1$. Denote by $H(z)$, for $|z|\leq 1$, be the probability generating function of $\{p_n\}_{n=1}^{\infty}$. Then, we have
\begin{equation}
\begin{split}
H(z)&=(e^{\alpha\delta}-1)^{-1}(\exp(-\alpha\log(1-\lambda z))-1)\\
\nonumber
&=\frac{c_0}{1-c_0}\sum_{n=1}^{\infty}{\alpha+n-1\choose\alpha-1}\lambda^nz^n.
\end{split}
\end{equation}
Thus, we obtain (4.3) for all $x \in \mathbb R_+.$\hfill $\Box$

\begin{theo} Let $0 < \lambda < 1$ and $\alpha >0$ and let $f(x)$ be a probability density function on $\mathbb R_+$. Define  a probability density function $p(x)$ by (4.1). Assume that $\int_0^{\infty}(f(x))^2dx < \infty.$ Then, the following are equivalent :

(a) $ p(x) \in \mathcal{S}_d$.

(b) $ f(x) \in \mathcal{S}_d$.

(c) $ f(x) \in \mathcal{L}_d$ and
$$ (1-c_0)p(x) \sim \varphi_{\alpha,\lambda}'(1)f(x).$$
\end{theo}
\medskip

Proof.  Assume that $\int_0^{\infty}(f(x))^2dx < \infty.$ Then, by (ii) and (iii) of Lemma 2.2, (b) implies (c). By Theorem 4.30 of Foss et al.  (2013), (c) implies (b) and, by (i) of Lemma 2.2, (b) implies  (a). Finally, we prove that (a) implies (b).   Define  a probability density functions $\phi(x)$ on $\mathbb R_+$  by (4.2). 
We find from Theorem 1.2 and Lemma 4.1 that $ \phi(x) \in \mathcal{S}_d$ and from (4.2) that  $\int_0^{\infty}(\phi(x))^2dx < \infty.$ 
We obtain from (ii) and (iii) of Lemma 2.2 that
\begin{equation}
(1-e^{-\alpha\delta})p(x) \sim \alpha\delta \phi(x).
\end{equation}
Define a function $f_0(x)$ on $\mathbb R_+$ as
$$f_0(x):=-\lambda^{-1}\sum_{n=1}^{\infty}\frac{(-\delta)^n}{n!}\phi^{n\otimes}(x).$$
Then, we have, for $t \in \mathbb R_+$,
\begin{equation}
\begin{split}
L_{f_0}(t)&=-\lambda^{-1}\sum_{n=1}^{\infty}\frac{(-\delta)^n}{n!}(L_{\phi}(t))^n\\ \nonumber
&=-\lambda^{-1}(\exp(-\delta L_{\phi}(t)) -1)\\
&=-\lambda^{-1}(\exp(\log(1-\lambda L_f(t)))-1)\\
&=L_f(t).
\end{split}
\end{equation}
Thus, we have $f_0(x)=f(x)$ for a.e. $x \in \mathbb R_+.$ Since
$\int_0^{\infty}(\phi(x))^2dx < \infty,$ we obtain from (ii) and (iii) of Lemma 2.2 that
\begin{equation}
f_0(x) \sim \lambda^{-1}\delta e^{-\delta}\phi(x).
\end{equation} 
Hence, by (i) of Lemma 2.2, we get
$$\tilde f_0(x):=f_0(x)\vee 0 \in \mathcal{S}_d.$$
Thus, we have by (4.4) and (4.5)
\begin{equation}
\begin{split}
c_0\alpha \lambda f(x)&=(1-c_0)p(x)-c_0\sum_{n=2}^{\infty}{\alpha+n-1\choose\alpha-1}\lambda^n \tilde f_0^{n\otimes}(x)\\ \nonumber 
&\sim \alpha\delta \phi(x)-(\varphi_{\alpha,\lambda}'(1)-c_0\alpha\lambda)\delta\lambda^{-1}(1-\lambda)\phi(x)\\
&=c_0\alpha\delta(1-\lambda)\phi(x).
\end{split}
\end{equation}
Hence, we see that
$$f(x) \sim \lambda^{-1}\delta (1-\lambda)\phi(x),$$
and from (i) of Lemma 2.2 that $f(x)\in \mathcal{S}_d.$\hfill $\Box$

\medskip

Next, we define the distribution $\pi$ of the supremum of a random walk. Let $\{X_n\}_{n=1}^{\infty}$ be IID random variables with common distribution $\rho$ on $\mathbb R$. Let $\{S_n\}_{n=0}^{\infty}$ be a random walk on $\mathbb R$ defined by $S_0:=0$ and $S_n:=\sum_{k=1}^{n}X_k$ for $n \geq 1$. Let $\pi$ be the distribution of the supremum $M$ of  $\{S_n\}$, that is, $M:=\sup_{n \geq 0} S_n.$ Define the measure $\nu$ on $(0,\infty)$ and a quantity $B$ as
\begin{equation}
\nu(dx):=1_{(0,\infty)}(x)\sum_{n=1}^{\infty}n^{-1}\rho^{n*}(dx)
\end{equation}
and
$$ B:=\sum_{n=1}^{\infty}n^{-1}P(S_n> 0)=\nu((0,\infty)).$$
It is well known that $M < \infty$ a.s. if and only if $B < \infty$ and that if  $B < \infty$, then $\pi$ is a compound Poisson distribution on $\mathbb R_+$ with L\'evy measure $\nu$. A sufficient condition for $M < \infty$ a.s. is that $-\infty < E(X_1) < 0.$ Define $\lambda$ as $\lambda:=1-e^{-B}$ when $B< \infty$. Let $Z^+$ be the first ascending ladder height in the random walk $\{S_n\}$ and denote the defective distribution on $\mathbb R_+$ by $\lambda \rho^+$. Then, $\rho^+$ is a distribution on $\mathbb R_+$ with $\rho^+(\{0\})=0.$ It is also well known that 
\begin{equation}
\pi= \sum_{n=0}^{\infty}(1-\lambda)\lambda^n(\rho^+)^{n*}.
\end{equation}
The distribution $\pi$ is important in classical ruin theory and queueing theory. See Asmussen and Albrecher (2010). Let $Z^-$ be the first descending ladder height in the random walk $\{S_n\}$ under the assumption that   $-\infty < E(X_1) < 0.$  We say that $\rho$ is {\it non-lattice} if the support of $\rho$ is not on any lattice. The following remark and the lemma are known up to now on the local subexponentiality of the distribution $\pi$.

\begin{rem} Suppose that  $-\infty < E(X_1) < 0$ and that $\rho$ is non-lattice.

(i) (Theorem 1 of Asmussen et al. (2002)) If $\rho \in \mathcal S^*,$ then we have
\begin{equation}
\pi((x,x+c]) \sim \frac{c}{|E(X_1)|}\bar\rho(x)
\end{equation}
for all $c >0$.

(ii) (Theorem 2 (b) of  Foss and Zachary (2003)) If $\rho \in\mathcal L$ and (4.8) holds for all $c >0$, then $\rho \in \mathcal S^*.$

(iii) (Lemma 3 of Asmussen et al. (2002)) If $\rho \in\mathcal L$, then
$$\rho^+((x,x+c]) \sim \frac{c}{\lambda |E(Z^-)|}\bar\rho(x)$$
for all $c >0$. Note that,
throughout in Asmussen et al. (2002),
$\rho \in \mathcal S^*$ is assumed but to prove Lemma 3 they did not use $\rho \in \mathcal S^*.$
\end{rem}
\begin{lem} (Theorem 6.2 of Watanabe and Yamamuro (2010)) Suppose that  $-\infty < E(X_1) < 0$ and that $\rho$ is non-lattice. Then, the following are equivalent : 

(1) $\pi \in \mathcal S_{loc}.$

(2)  $\nu_{(1)} \in \mathcal S_{loc}.$

(3)  $\rho^+ \in \mathcal S_{loc}.$

(4) $\rho \in \mathcal S^*.$
 
 (5) $\rho \in\mathcal L$ and (4.8) holds.
 
 (6)  $\rho \in\mathcal L$ and there is $C \in (0,\infty)$ such that
 $\pi((x,x+c]) \sim Cc\bar\rho(x)$
 for all $c > 0$.
\end{lem}

We find from (4.6) and (4.7) that  $\rho$ is absolutely continuous on $(0,\infty)$ if and only if  so is $\rho^+$ on $\mathbb R_+$. Let  $\rho^+(dx):=f(x)dx$ and
$\lambda:=1-e^{-B}\in (0,1).$ Define
$p(x)$ and $\phi(x)$ on $\mathbb R_+$ as
\begin{equation}
 p(x):=(1-\lambda)\sum_{n=1}^{\infty}\lambda^{n-1}f^{n\otimes}(x)
 \end{equation}
and
 \begin{equation}
 \phi(x):=B^{-1}\sum_{n=1}^{\infty}\frac{\lambda^n}{n}f^{n\otimes}(x).
 \end{equation}
Then, we have by Lemma 4.1 with $\alpha=1$ and $\delta=B$
$$ p(x)=(e^{B}-1)^{-1}\sum_{n=1}^{\infty}\frac{B^n}{n!}\phi^{n\otimes}(x)$$
for $x \in \mathbb R_+$. Moreover, $\pi$ and $\nu_{(0)}$ are represented on $\mathbb R_+$ as
$$ \pi(dx)=e^{-B}\delta_0(dx)+(1-e^{-B})p(x)dx$$
and $\nu_{(0)}(dx)=\phi(x)dx.$

\begin{theo} Suppose that  $-\infty < E(X_1) < 0$. Assume that $\rho$ is absolutely continuous on $(0,\infty)$. Let $p(x)$ and $\phi(x)$ be given by (4.9) and (4.10), respectively. Further, assume that $\int_0^{\infty}(p(x))^2dx < \infty.$ Then, the following are equivalent : 

(1) $p(x) \in \mathcal S_{d}.$

(2)  $\phi(x) \in \mathcal S_{d}.$

(3)  $f(x) \in \mathcal S_{d}.$

(4) $\rho \in \mathcal S^*$ and $\phi(x) \in \mathcal L_{d}.$

 (5) $\rho \in\mathcal L$ and 
 $$(1-e^{-B})p(x) \sim \frac{\bar\rho(x)}{|E(X_1)|}.$$
  
 (6)  $\rho \in\mathcal L$ and there is $C \in (0,\infty)$ such that
 $(1-e^{-B})p(x) \sim C\bar\rho(x).$
 \end{theo}

\medskip

Proof. Note from (4.9) and (4.10) that the three conditions $\int_0^{\infty}(p(x))^2dx < \infty,$ $\int_0^{\infty}(\phi (x))^2dx < \infty.$ and $\int_0^{\infty}(f(x))^2dx < \infty$ are equivalent. Equivalence of (1) and (2) is due to Theorem 1.1 and that of (1) and (3) is due to Theorem 4.1. Suppose that (4) holds. Then, we find from Lemma 4.2 that $\phi(x)dx \in \mathcal S_{loc}$ and hence, by $\phi(x) \in \mathcal L_{d}$ and (ii) of Lemma 2.1, (2) holds. Conversely, by Lemma 4.2, (2) implies (4). Suppose that (1) holds. Then, we obtain from Lemma 4.2 that $\rho \in\mathcal L$ and 
 $$(1-e^{-B})p(x)\sim \pi((x,x+1]) \sim \frac{\bar\rho(x)}{|E(X_1)|}.$$
   Thus, (5) and (6) hold. Conversely, suppose that (6) holds. Then, 
   we have $p(x) \in \mathcal L_{d}.$ We see from Lemma 4.2 that $C=1/|E(X_1)|$ and (5) holds and $\pi \in \mathcal S_{loc}.$ Thus, by (ii) of Lemma 2.1, $p(x) \in \mathcal S_{d}.$
Therefore, all assertions (1)-(6) are equivalent.\hfill $\Box$
\begin{rem} Our result characterizes the absolutely continuous case. On the other hand, Korshunov (2006) discussed the non-absolutely continuous case and showed in his  Theorem 4 that if $\rho \in \mathcal S^*$ and the absolutely continuous part of $\rho$ is in  $\mathcal L_{d}$, then assertion (5) of the above theorem holds for the absolutely continuous part of the distribution  $\pi$ under an additional assumption for the singular part of $\rho$. However, he did not obtain the converse direction.
\end{rem}


\medskip
\begin{thebibliography}{99}
\bibitem{aa} Asmussen, S., Albrecher, H.: Ruin probabilities. Second edition. Advanced Series on Statistical Science and Applied Probability, 14 World Scientific Publishing Co. Pte. Ltd., Hackensack, NJ,(2010).
\bibitem{afk} Asmussen, S., Foss, S., Korshunov, D.: Asymptotics for sums of random variables with local subexponential behaviour, J. Theoret. Probab. 16,  489-518 (2003). 
\bibitem{akkkt} Asmussen, S., Kalashnikov, V., Konstantinides, D., Kl\"uppelberg, C., Tsitsiashvili, G.: A local limit theorem for random walk maxima with heavy tails. Statist. Probab. Lett. 56, 399-404 (2002).
\bibitem{cnw} Chover, J., Ney, P., Wainger, S.: Functions of probability measures. J. Analyse Math. 26, 255-302 (1973). 
\bibitem{ft} Finkelshtein, D.,  Tkachov, P.: Kesten's bound for subexponential densities on the real line and its multi-dimensional analogues. Adv.  Appl. Probab. 50, 373-395 (2018).
\bibitem{fkz} Foss, S., Korshunov, D., Zachary, S.: An introduction to heavy-tailed and subexponential distributions. Second edition. Springer Series in Operations Research and Financial Engineering. Springer, New York,(2013).
\bibitem{fz} Foss, S., Zachary, S.: The maximum on a random time interval of a random walk with long-tailed increments and negative drift. Ann. Appl. Probab. 13, 37-53 (2003). 
\bibitem{k1} Kl\"uppelberg, C.: Subexponential distributions and integrated tails. J. Appl. Probab. 25, 132-141 (1988). 
\bibitem{k2}  Kl\"uppelberg, C.: Subexponential distributions and characterizations of related classes. Probab. Theory Related Fields 82, 259-269 (1989). 
 \bibitem{ko} Korshunov, D.: On the distribution density of the supremum of a random walk in the subexponential case. Siberian Math. J.
 47, 1060–1065 (2006).
 \bibitem{s} Sato, K.: L\'evy processes and infinitely divisible distributions. Cambridge Studies in Advanced Mathematics, 68 Cambridge Univ. Press.(2013).
\bibitem{w} Watanabe, T.: Convolution equivalence and distributions of random sums. Probab. Theory Related Fields 142, 367-397 (2008). 
 
\bibitem{wy2} Watanabe, T., Yamamuro, K.: Local subexponentiality and \\
self-decomposability. J. Theoret. Probab. 23, 1039-1067 (2010). 
\bibitem{wy4} Watanabe, T., Yamamuro, K.: Two non-closure properties on the class of subexponential densities. J. Theoret. Probab. 30, 1059-1075 (2017).
\end{thebibliography}
\end{document}